\documentclass[12pt]{article}
\usepackage{latexsym, amssymb, amscd, amsmath}
\begin{document}
\title{\textbf{The Geometry of the Moduli Space \newline of Curves of 
Genus $23$ }}
\author{GAVRIL FARKAS
}
\date{ }
\maketitle

\section{Introduction}
The problem of describing the birational geometry of the moduli  space
$\mathcal{M}_g$ of complex curves of genus $g$ has a long history. Severi
already knew in 1915 that $\mathcal{M} _g$ is unirational for $g\leq 10$
(cf. \cite{Sev}; see also \cite{AC1} for a modern proof). In the same
paper Severi conjectured that $\mathcal{M} _g$ is unirational for all
genera $g$. Then for a long period this problem seemed intractable
(Mumford writes in \cite{Mu}, p.51:\lq\lq Whether more $\mathcal{M}_g$'s,
$g\geq 11$, are unirational or not is a very interesting problem, but one
which looks very hard too, especially if $g$ is quite large"). The
breakthrough came in the eighties when Eisenbud, Harris and Mumford proved
that $\mathcal{M} _g$ is of general type as soon as $g\geq 24$ and that
the Kodaira dimension of $\mathcal{M} _{23}$ is $\geq 1$ (see \cite{HM},
\cite{EH3}). We note that $\mathcal{M}_g$ is rational for $g\leq 6$ (see
\cite{Dol} for problems concerning the rationality of various moduli
spaces).
\newline
\indent Severi's proof of the unirationality of $\mathcal{M}_g$ for small 
$g$ was based on representing a general curve of genus $g$ as a plane
curve of degree $d$ with $\delta $ nodes; this is possible when $d\geq
2g/3+2.$ When the number of nodes is small, i.e. $\delta \leq
(d+1)(d+2)/6$, the dominant map from the variety of plane curves of degree
$d$ and genus $g$ to $\mathcal{M}_g$ yields a rational parametrization of
the moduli space. The two conditions involving $d$ and $\delta $ can be
satisfied only when $g\leq 10$, so Severi's argument cannot be extended
for other genera. However, using much more subtle ideas, Chang, Ran and
Sernesi proved the unirationality of $\mathcal{M}_g$ for $g=11,12,13$ (see
\cite{CR1}, \cite{Se1}), while for $g=15,16$ they proved that the Kodaira
dimension is $-\infty $ (see [CR2,4] ). The remaining cases $g=14$ and
$17\leq g\leq 23$ are still quite mysterious. Harris and Morrison
conjectured in \cite{HMo} that $\mathcal{M}_g$ is uniruled precisely when
$g<23$.
\newline
\indent All these facts indicate that $\mathcal{M}_{23}$ is a very 
interesting transition case. Our main result is the following:
\newtheorem{thm}{Theorem}
\begin{thm}
The Kodaira dimension of the moduli space of curves of genus $23$ is 
$\geq 2$.
\end{thm}
We will also present some evidence for the hypothesis that the Kodaira 
dimension of $\mathcal{M}_{23}$ is actually equal to $2$.\newline
\textbf{Acknowledgments: } I am grateful to my advisor Gerard van der Geer, for proposing the problem and for continuous support and guidance. I thank Joe Harris for many inspiring discussions during my stay at Harvard. I also benefited from conversations with Marc Coppens and Carel Faber.
\section{Multicanonical linear systems and the Kodaira dimension of 
$\mathcal{M}_g$} We study three multicanonical divisors on
$\mathcal{M}_{23}$, which are (modulo some boundary components) of
Brill-Noether type, and we conclude by looking at their relative position
that $\kappa (\mathcal{M}_{23})\geq 2$.\newline
\indent We review some notations. We shall denote by 
$\overline{\mathcal{M}}_g$ and $\overline{\mathcal{C}}_g$ the moduli
spaces of stable and $1$-pointed stable curves of genus $g$ over $\mathbb
C$. If $C$ is a smooth algebraic curve of genus $g$, we consider for any
$r$ and $d$, the scheme whose points are the $\mathfrak g^r_d$'s on $C$,
that is,
$$G^r_d(C)=\{ (\mathcal{L},V):\mathcal{L} \in \mbox{Pic }^d(C), V\subseteq
H^0(C,\mathcal{L} ),\mbox{dim}(V)=r+1\},$$  (cf. \cite{ACGH}) and denote
the associated Brill-Noether locus in $\mathcal{M} _g$ by
$$\mathcal{M} ^r_{g,d} :=\{ [C]\in \mathcal{M} _g:G^r_d(C)\neq 
\emptyset\},$$ and by $\overline{\mathcal{M} } _{g,d}^r$ its closure in
$\overline{\mathcal{M} } _g$.\newline
\indent The distribution of linear series on algebraic curves is  governed
(to some extent) by the {\sl Brill-Noether number}
$$\rho (g,r,d):=g-(r+1)(g-d+r).$$
The Brill-Noether Theorem asserts that when $\rho (g,r,d)\geq 0$  every
curve of genus $g$ possesses a $\mathfrak g ^r_d$, while when $\rho
(g,r,d)<0$ the general curve of genus $g$ has no $\mathfrak g ^r_d$'s,
hence in this case the Brill-Noether loci are proper subvarieties of
$\mathcal{M}_g.$ When $\rho (g,r,d)<0$, the naive expectation that $-\rho
(g,r,d)$ is the  codimension of $\mathcal{M}_{g,d}^r$ inside
$\mathcal{M}_g$, is in general way off the mark, since there are plenty of
examples of Brill-Noether loci of unexpected dimension (cf. \cite{EH2}).
However, we have Steffen's result in one direction (see \cite{St}):
\newline
\newline
\textit{ If $\rho (g,r,d)<0$ then each component of $\mathcal{M}_{g,d}^r$
has  codimension at most $-\rho (g,r,d)$ in ~$\mathcal{M}_g.$}
\newline
\newline
On the other hand, when the Brill-Noether number is not very negative, 
the Brill-Noether loci tend to behave nicely. Existence of components of
$\mathcal{M}_{g,d}^r$ of the expected dimension has been proved for a
rather wide range, namely for those $g, r, d$ such that $\rho (g,r,d)<0$,
and
$$\rho (g,r,d)\geq  
\begin{cases}
-g+r+3  &\mbox {  if } r\mbox{ is odd; } \\
-rg/(r+2)+r+3  &\mbox{  if } r\mbox{ is even.} 
\end{cases}
$$ 
We have a complete answer only when $\rho (g,r,d)=-1$. Eisenbud and Harris have
proved in \cite{EH2} that in this case $\mathcal{M}_{g,d}^r$ has a unique 
divisorial component, and using the previously mentioned theorem of
Steffen's, we obtain the following result:
\begin{center}
\textit{ If $\rho (g,r,d)=-1$, then $\overline{\mathcal{M}}_{g,d}^r$ is
an  irreducible divisor of $\overline{\mathcal{M}}_g$.}
\end{center}
We will also need Edidin's result (see \cite{Ed2} ) which says that for $g\geq 12$ and $\rho(g, r, d)=-2$, 
all components of $\mathcal{M}_{g,d}^r$ have codimension $2$. We can get
codimension $1$ Brill-Noether conditions only for the genera $g$ for which
$g+1$ is composite. In that case we can write
$$
g+1=(r+1)(s-1),\mbox{ }s\geq 3
$$
and set $d:=rs-1.$ Obviously $\rho (g,r,d)=-1$ and 
$\overline{\mathcal{M}}_{g,d}^r$ is an irreducible divisor. Furthermore,
its class has been computed (cf. \cite{EH3} ):
$$
[\overline{\mathcal{M}}_{g,d}^r]=c_{g,r,d}
\left((g+3)\lambda-\frac{g+1}{6}\delta
_{0}-\sum_{i=1}^{[g/2]} i(g-i)\delta _{i}\right),
$$ 
where $c_{g,r,d}$ is a
positive rational number equal to $3\mu /(2g-4)$, with $\mu $ being the
number of $\mathfrak g ^{r}_{d}$'s on a general pointed curve $(C_0,q)$ of
genus $g-2$ with ramification sequence $(0,1,2,\ldots ,2)$ at $q$.  For
$g=23$ we have the following possibilities:
$$
(r,s,d)=(1,13,12),\mbox{ } (11,3,32),\mbox{ }(2,9,17),\mbox{ }(7,4,24),\mbox{
}(3,7,20),\mbox{ }(5,5,24).
$$ 
It is immediate by Serre duality, that cases
$(1,13,12)$ and $(11,3,32)$  yield the same divisor on 
$\mathcal{M}  _{23}$, namely the $12$-gonal locus $\mathcal{M} _{12}^1$;
similarly, cases $(2,9,17)$ and $(7,4,24)$ yield the divisor $\mathcal{M}
^{2}_{17}$ of curves having a $\mathfrak g ^2_{17}$, while cases
$(3,7,20)$ and $(5,5,24)$ give rise to $\mathcal{M} ^{3}_{20}$, the
divisor of curves having a $\mathfrak g ^3_{20}$. Note that when the genus
we are referring to is clear from the context, we write $\mathcal{M}
^r_d=\mathcal{M} _{g,d}^r.$
\newline
\indent By comparing the classes of the Brill-Noether divisors to the class of 
the canonical divisor
$ K_{\overline{\mathcal{M} }_{g,reg}}=13\lambda -2\delta _0 -3\delta _1
-\cdots -2\delta _{[g/2]},$
at least in the case when $g+1$ is composite we can infer that
$$K_{\overline{\mathcal{M} }_{g,reg}}=a[\overline{\mathcal{M} }_{g,d}^r]+b\lambda
+(\mbox{ positive combination of } \delta _0,\ldots ,\delta _{[g/2 ]} ),$$
where $a$ is a positive rational number, while $b>0$ as long as $g\geq 24$ 
but $b=0$ for $g=23.$
As it is well-known that $\lambda $ is big on $\overline{\mathcal{M}} _g$,  it
follows that $\mathcal{M}_g$ is of general type for $g\geq 24$ and that it has
non-negative Kodaira dimension when $g=23.$ Specifically for $g=23$, we get that
there are positive integer  constants $m, m_1, m_2, m_3$ such that:
\begin{equation}
 mK=m_1[\overline{\mathcal{M}}^{1}_{12}]+E,\mbox{ }  mK=m_2[\overline{\mathcal{M}
}^{2}_{17}]+E,\mbox{ } mK=m_3[\overline{\mathcal{M}}^{3}_{20}]+E,
\end{equation}
where $E$ is the same positive combination of $\delta _1,\ldots ,\delta _{11}$.
\newtheorem{prop}{Proposition}[section]
\begin{prop}[Eisenbud-Harris, \cite{EH3}]
 There exists a smooth curve of genus $23$ that possesses a $\mathfrak g^1_{12}$ , but
no $\mathfrak g^2_{17}$. It follows that $\kappa  (\mathcal{M} _{23})\geq 1.$
\end{prop}
 \indent Harris and Mumford proved (cf. \cite{HM}) that $\overline{\mathcal{M} }_g$
has only canonical singularities for $g\geq 4$, hence 
$ H^0(\overline{\mathcal{M}}_{g,reg},nK)=H^0 (\widetilde{ \mathcal{M}}_{g},nK)$ for
each $n\geq 0,$ with
$\widetilde{\mathcal{M}_g}$ a desingularization of $\overline{\mathcal{M}}_g$. We
already know that $\mbox{dim}(\mbox{Im}\phi_{mK})\geq 1$, where
 $\phi _{mK}:\overline{\mathcal{M} }_{23}\displaystyle \mathop {-- \rightarrow}
\mathbb{P}^{\nu }$ is the multicanonical map, $m$ being as in (1). We will prove that
$\kappa (\mathcal{M}_{23})\geq 2$. Indeed, let us assume that
$\mbox{dim}(\mbox{Im}\phi_{mK})=1$. Denote by $C:=\overline{\mbox{Im}\phi_{mK}}$ the
Kodaira image of $\overline{\mathcal{M} }_{23}$. We reach a contradiction by proving
two things:
\newline
$\bullet$ $\alpha )$ The Brill-Noether divisors $\mathcal{M}^1_{12},
\mathcal{M}^2_{17}$ and $\mathcal{M}^3_{20}$ are mutually distinct.\newline 
$\bullet$ $\beta)$
There exist smooth curves of genus $23$ which belong to exactly two of the 
Brill-Noether divisors from above.
\newline
This suffices in order to prove Theorem 1: since $\overline{\mathcal{M}}^1_{12}$,
$\overline{\mathcal{M} }^2_{17}$ and $\overline{\mathcal{M}}^3_{20}$ are part of
different multicanonical divisors, they must be contained in different fibres of the
multicanonical map $\phi_{mK}$. Hence there exists different points $x,y,z\in C$ such
that $$ \mathcal{M}^1_{12}=\overline{\phi ^{-1}(x)} \cap \mathcal{M} _{23}, \mbox{ }
\mathcal{M}^2_{17} =\overline{\phi^{-1}(y)} \cap \mathcal{M}_{23},\mbox{ }
\mathcal{M}^3_{20} =\overline{\phi^{-1}(z)} \cap \mathcal{M}_{23}.$$ It follows that
the set-theoretic intersection of any two of them will be contained in the base locus
of $|mK_{\overline{\mathcal{M} }_{23}}|$. In particular:
\begin{equation}
\mbox{supp}(\mathcal{M}^1_{12})\cap \mbox{supp}(\mathcal{M} ^2
_{17})=\mbox{supp}(\mathcal{M}^2_{17})\cap \mbox{supp}(\mathcal{M}^3_{20})=
\mbox{supp}(\mathcal{M} ^3_{20})\cap \mbox{supp}(\mathcal{M}^1_{12}),
\end{equation}
and this contradicts $\beta )$. We complete the proof of $\alpha) $ and $\beta) $ 
is Section 5.
%%
%%%
%%%
%%
\section{Deformation theory for $\mathfrak g ^r_d$'s and limit linear series}
We recall a few things about the variety parametrising $\mathfrak g^r_d$'s on the 
fibres of the universal curve (cf. \cite{AC2}), and then we recap on the theory of
limit linear series (cf. \cite{EH1}, \cite{Mod}), which is our main technique for the
study of $\mathcal{M} _{23}$.
\newline \indent Given $g, r, d$ and a point $[C]\in \mathcal{M}_g$, there is a 
connected neighbourhood $U$ of $[C]$, a finite ramified covering 
$h:\mathcal{M}\rightarrow U$, such that $\mathcal{M}$ is a fine moduli space of curves
(i.e. there exists $\xi :\mathcal{C}\rightarrow \mathcal{M}$ a universal curve), and a
proper variety over $\mathcal{M}$,
$$\pi :\mathcal{G} ^r_d\rightarrow \mathcal{M} $$
which parametrizes classes of couples $(C,l)$, with $[C]\in \mathcal{M}$ and $l\in G^r_d(C)$, where we have made the identification $C=\xi ^{-1}([C])$ .\newline
Let $(C,l)$ be a point of $\mathcal{G} ^r_d$ corresponding to a curve $C$ and a
linear series $l=(\mathcal{L} ,V),$ where $\mathcal{L} \in \mbox{Pic}^d(C),V\subseteq H^0(C, \mathcal{L} ), \mbox{ and dim}(V)=r+1.$ 
By choosing a basis in $V$, one has a morphism $f:C\rightarrow \mathbb P ^r.$ The normal sheaf of $f$ is defined through the exact sequence
\begin{equation}
0\longrightarrow T_C\longrightarrow f^{*}(T_{\mathbb P ^r})\longrightarrow
N_f\longrightarrow 0.
\end{equation}
By dividing out the torsion of $N_f$ one gets to the exact sequence 
\begin{equation}
0\longrightarrow \mathcal{K} _f\longrightarrow N_f\longrightarrow N_f'
\longrightarrow 0,
\end{equation}
where the torsion sheaf $\mathcal{K} _f$ (the cuspidal sheaf) is based at those points $x\in C$ where $df(x)=0$, and $N_f'$ is locally free of rank $r-1.$
The tangent space $T_{(C, l)}(\mathcal{G} ^r_d)$ fits into an exact
sequence (cf. \cite{AC2}):
\begin{equation}
0\longrightarrow \mathbb C \longrightarrow \mbox{Hom}(V,V)\longrightarrow H^0(C,
N_f)\longrightarrow T_{(C, l)}(\mathcal{G} ^r_d)\longrightarrow 0,
\end{equation}
from which we have that $\mbox{dim }T_{(C,l)}(\mathcal{G} ^r_d)=3g-3+\rho(g,r,d)+h^1(C,N_f).$
\begin{prop}
Let $C$ be a curve and $l\in G^r_d(C)$ a base point free linear series. Then the variety $\mathcal{G} ^r_d$ is smooth and of dimension $3g-3+\rho (g,r,d)$ at the point $(C, l)$ if and only if $H^1(C, N_f)=0.$
\end{prop}
\textbf{Remark: }
 The condition $H^1(C, N_f)=0$ is automatically satisfied for $r=1$
as $N_f$ is a sheaf with finite support. Thus $\mathcal{G}^ 1_d$ is smooth
of dimension $2g+2d-5.$ It follows that $\mathcal{G} ^1_d$ is birationally
equivalent to the $d$-gonal locus $\mathcal{M}^1_d$ when $d<(g+2)/2 .$

\vskip10pt

Limit linear series try to answer questions of the following kind: what happens to a family of $\mathfrak g ^r_d$'s when a smooth curve specializes to a reducible curve? Limit linear series solve such problems for a class of reducible curves, those of compact type. A curve $C$ is of {\sl compact type } if the dual graph is a tree. A curve $C$ is {\sl tree-like } if, after deleting edges leading from a node to itself, the dual graph becomes a tree.
\newline
\indent Let $C$ be a smooth curve of genus $g$ and $l=(\mathcal{L} ,V)\in G^r_d(C)$, $\mathcal{L} \in \mbox{Pic}^d(C), V\subseteq H^0(C,\mathcal{L} ), \mbox{ and dim}(V)=r+1.$
Fix $p\in C$ a point. By ordering the finite set $\{ord_p(\sigma ) \}_{\sigma \in V}$
one gets the {\sl vanishing sequence} of $l$ at $p$:
$$a^l(p):0\leq a_0^l(p)<\ldots <a_r^l(p)\leq d.$$
The {\sl ramification sequence} of $l$ at $p$ 
$$\alpha ^l(p):0\leq \alpha ^l_{0}(p)\leq \ldots \leq \alpha ^l_{r}(p)\leq d-r$$
is defined as
$\alpha ^l_{i}(p)=a^l_{i}(p)-i$ and the {\sl weight} of $l$ at $p$ is 
$$w^l(p)=\sum_{i=0}^{r} \alpha _i^l(p) .$$
A {\sl Schubert index of type $(r, d)$ } is a sequence of integers
$\beta :0\leq \beta _0\leq \ldots \beta _r\leq d-r.$
If $\alpha $ and $\beta $ are Schubert indices of type $(r, d)$ we write
$\alpha \leq \beta \Longleftrightarrow \alpha _i\leq \beta _i, i=0,\ldots ,r.$
The point $p$ is said to be a {\sl ramification point } of $l$ if $w^l(p)>0.$
The linear series $l$ is said to have a {\sl cusp } at $p$ if $\alpha ^l(p)\geq (0,1,
\ldots ,1).$
For $C$ a tree-like curve, $p_1,\ldots ,p_n\in C$ smooth points and 
$\alpha ^1,\ldots ,\alpha ^n$ Schubert indices of type $(r,d)$, we define
$$G^r_d(C, (p_1, \alpha ^1),\ldots (p_n, \alpha ^n)):=\{ l\in G^r_d(C):
\alpha ^l(p_1)\geq \alpha ^1,\ldots ,\alpha ^l(p_n)\geq \alpha ^n\}.$$
This scheme can be realized naturally as a determinantal variety and its 
expected dimension is 
$$\rho (g, r, d, \alpha ^1,\ldots ,\alpha ^n):=\rho(g, r, d)-\sum_{i=1}^{n} 
\sum_{j=0}^{r} \alpha ^i_j.$$
If $C$ is a curve of compact type, a {\sl crude limit $\mathfrak g ^r_d$ } on $C$ is a
collection of ordinary linear series
$l=\{l_Y\in G^r_d(Y): Y\subseteq C\mbox{ is a component } \}$,
satisfying the following compatibility condition:
if $Y$ and $Z$ are components of $C$, with $\{ p \} =Y\cap Z $, then
$$a^{l_Y}_{i}(p)+a^{l_Z}_{r-i}(p)\geq d,\mbox{ for }i=0,\ldots r.$$
If equality holds everywhere, we say that $l$ is a
{\sl refined limit $\mathfrak g ^r_d$.} The `honest' linear series $l_Y\in G^r_d(Y)$ is
called the $Y$-aspect of the limit linear series $l$.
\newline
\indent We will often use the additivity of the Brill-Noether number: if $C$ is a curve of compact type,
for each component $Y\subseteq C$, let $q_1,\ldots ,q_s$ be the points where $Y$ meets the
other components of $C$. Then for any limit
$\mathfrak g ^r_d$ on $C$ we have the following inequality:
\begin{equation}
\rho (g,r,d)\geq \sum_{Y\subseteq C} \rho (l_Y, \alpha ^{l_Y}(q_1),\ldots ,
\alpha ^{l_Y}(q_s)),
\end{equation}
with equality if and only if $l$ is a refined limit linear series.
\newline
\indent It has been proved in \cite{EH1} that limit linear series arise indeed as limits of
ordinary linear series on smooth curves. Suppose we are given a family $\pi :\mathcal{C}
\rightarrow B$ of genus $g$ curves, where 
$B= \mbox{Spec} (R)$ with $R$ a complete discrete valuation ring.  
 Assume furthermore that $\mathcal{C} $ is a smooth surface and that if 
$0, \eta $ denote the special and generic point of $B$ respectively, 
the central fibre $C_0$ is reduced and of compact type, while the generic geometric fibre
$C_{\eta }$ is smooth and irreducible. If $l_{\eta }=(\mathcal{L} _{\eta }, V_{\eta })$ is a
$\mathfrak g ^r_d$ on
$X_{\eta}$, there is a canonical way to associate a crude limit series $l_0$ on $C_0$ which is
the limit of $l_{\eta }$ in a natural way: for each component $Y$ of $C_0$, there exists a
unique line bundle $\mathcal{L} ^Y$ on $\mathcal{C} $ such that $$\mathcal{L}^Y_{|\mathcal{C}
_{\eta }}=\mathcal{L}_{\eta }\mbox{ and }\mbox{ }  \mbox{deg}_Z(\mathcal{L}^Y_{|_{Z}})=0,$$
for any component $Z$ of $C_0$ with
$Z\neq Y$. (This implies of course that $\mbox{deg}_Y(\mathcal{L} ^Y_{|_Y})=d).$ Define
$V^Y=V_{\eta }\cap H^0(\mathcal{C} ,\mathcal{L}^Y)\subseteq H^0(\mathcal{C} _{\eta
},\mathcal{L} _{\eta }).$ Clearly, $V^Y$ is a free $R$-module of rank $r+1.$
\newline
Moreover, the  composite homomorphism 
$$
V^Y(0)\rightarrow (\pi _{*}\mathcal{L} ^Y)(0)\rightarrow H^0(C_0,\mathcal{L}
 ^Y_{|_{C_0}})\rightarrow H^0(Y, \mathcal{L} ^Y_{|_{C_0}})
$$
is injective, hence $l_Y=(\mathcal{L} ^Y_{|_Y}, V^Y(0))$ is an ordinary $\mathfrak g ^r_d$
on $Y$. 
One proves that
$l=\{ l_Y:Y \mbox{ component of }C_0 \}$ is a limit linear series.
\newline  
 \indent If $C$ is a reducible curve of compact type, 
$l$ a limit $\mathfrak g ^r_d$ on $C$, we say that  $l$ is {\sl smoothable } if there
exists $\pi :\mathcal{C} \rightarrow B$ a family of curves with central fibre $C=C_0$
as above, and $(\mathcal{L} _{\eta },V_{\eta })$ a $\mathfrak g^r_d$ on the generic
fibre $C_{\eta }$ whose limit on $C$ (in the sense previously described) is $l$.
\newline
\textbf{Remark: } If a stable curve of compact type $C$, has no limit $\mathfrak g ^r_d$'s, 
then $[C] \notin \overline{\mathcal{M} } _{g,d}^r.$  If there exists a smoothable limit
$\mathfrak g ^r_d$ on $C$, then $[C]\in \overline{\mathcal{M} } _{g,d}^r.$
\newline
\indent Now we explain a criterion due to Eisenbud and Harris 
(cf. \cite{EH1}), which gives a sufficient condition for a limit $\mathfrak g^r_d$ to be 
smoothable. Let $l$ be a limit $\mathfrak g ^r_d$ on a curve $C$ of compact type. Fix
$Y\subseteq C$ a component, and $\{q_1, \ldots ,q_s\}=Y\cap \overline{(C-Y)}.$ Let
$$\pi :\mathcal{Y} \rightarrow B,\mbox{ } \tilde {q_i}:B\rightarrow \mathcal{Y}$$
be the versal deformation space of $(Y, q_1,\ldots q_s).$ The base $B$ can be viewed as  a small
$(3g(Y)-3+s)$-dimensional polydisk. Using general theory one constructs a proper scheme over $B$,
$$\sigma :\mathcal{G} ^r_{d}(\mathcal{Y} /B; (\tilde{q} _i,\alpha ^{l_Y}(q_i))_{i=1}^s)
\rightarrow B$$ whose fibre over each $b\in B$ is
$\sigma ^{-1}(b)=G^r_d(Y_b, (\tilde{q} _i(b),\alpha ^{l_Y}(q_i))_{i=1}^s).$
One says that $l$ is {\sl dimensionally proper with respect to $Y$}, if the
$Y$-aspect $l_Y$ is contained in some component $\mathcal{G} $ of 
$\mathcal{G} ^r_d(\mathcal{Y} /B; (\tilde{q} _i, \alpha ^{l_Y}(q_i))_{i=1}^s)$  of the expected
dimension, i.e.
$$\mbox{dim } \mathcal{G} =\mbox{dim } B +\rho(l_Y, \alpha ^{l_Y}(q_1), \ldots
\alpha^{l_Y}(q_s)).$$ One says that $l$ is {\sl dimensionally proper}, if it is dimensionally
proper with respect to any component $Y\subseteq C.$ The `Regeneration Theorem' (cf. \cite{EH1})
states that every dimensionally proper limit linear series is smoothable.
\newline
\indent The next result is a `strong Brill-Noether Theorem', i.e. it not  only asserts a
Brill-Noether type statement, but also singles out the locus where the statement fails. 
\begin{prop}[Eisenbud-Harris]
Let $C$ be a tree-like curve and for any component $Y\subseteq C$, denote
by $q_1, \ldots ,q_s\in Y$ the points where $Y$ meets the other components of
$C$. Assume that for each $Y$ the following conditions are satisfied:
\begin{itemize}
\item[a.] If $g(Y)=1$ then $s=1.$
\item[b.] If $g(Y)=2$ then $s=1$ and $q$ is not a Weierstrass point.
\item[c.] If $g(Y)\geq 3$ then $(Y, q_1, \ldots ,q_s)$ is a general $s$-pointed curve.
\end{itemize}
Then for $p_1, \ldots p_t\in C$ general points, $\rho (l,\alpha ^l(p_1), \ldots , \alpha
^l(p_t))\geq 0$ for any limit linear series on $C$.
\end{prop}
Simple examples involving pointed elliptic curves show that the 
condition $\rho (g,r,d)\geq \sum_{i=1}^{t} w^l(p_i)$ does not guarantee the existence  of a
linear series $l\in G^r_d(C)$ with prescribed ramification at general points $p_1,p_2,\ldots
,p_t\in C$. The appropriate condition in the pointed case can be given in terms of Schubert
cycles. Let $\alpha =(\alpha _0,\ldots ,\alpha _r)$ be a Schubert index of type $(r, d)$ and 
$$\mathbb C^{d+1}=W_0\supset W_1\supset \ldots \supset W_{d+1}=0$$ a decreasing flag of  linear
subspaces. We consider the Schubert cycle in the Grassmanian,
$$\sigma _{\alpha }=\{ \Lambda \in G(r+1, d+1): \mbox{dim} (\Lambda \cap W_{\alpha _i+
i})\geq r+1-i,\mbox{ } i=0,\ldots ,r \} .$$
For a general $t$-pointed curve $(C,p_1,\ldots ,p_t)$ of genus $g$, and  
$\alpha ^1,\ldots ,\alpha ^t$ Schubert indices of type $(r,d)$, the
necessary and sufficient condition that $C$ has a $\mathfrak g ^r_d$ with ramification
$\alpha ^i$ at $p_i$ is that
\begin{equation}
\sigma _{_{\alpha ^1}}\cdot \ldots \cdot \sigma _{_{\alpha ^t}}\cdot \sigma_{_{(0,1,\ldots ,1)}}^g\neq 0 \mbox{ in } H^*(G(r+1,d+1),\mathbb Z) .
\end{equation}
In the case $t=1$ this condition can be made more explicit (cf. \cite{EH3}): a general pointed curve $(C,p)$ of genus $g$ carries a $\mathfrak g ^r_d$ with
ramification sequence $(\alpha _0,\ldots ,\alpha _r)$ at $p$, if and only if
\begin{equation} 
\sum_{i=0}^r (\alpha _i +g-d+r)_+\leq g,
\end{equation}
where $x_+=\mbox{max}\{x,0\}.$ One can make the following simple but useful observation:
\begin{prop}
Let $(C, p, q )$ be a general $2$-pointed curve of genus $g$ and $(\alpha _0,
\ldots , \alpha _r)$ a Schubert index of type $(r, d)$. Then $C$ has a $\mathfrak g ^r_d$
having ramification sequence $(\alpha _0, \ldots ,\alpha _r)$ at $p$ and a cusp
at $q$ if and only if 
$$\sum_{i=0}^r (\alpha _i+g+1-d+r)_+\leq g+1.$$
\end{prop}
\textsl{Proof: } The condition for the existence of the $\mathfrak g ^r_d$ with
ramification $\alpha $ at $p$ and a cusp at $q$ is that
$\sigma _{\alpha }\cdot \sigma _{_{(0,1,\ldots ,1)}}^{g+1}\neq 0$ (cf. (7)).
According to the Littlewood-Richardson rule (see \cite{F}), this is equivalent
with 
$\sum_{i=0}^r (\alpha _i+g+1-d+r)_+\leq g+1.\mbox{  } \hfill \Box$
\section{A few consequences of limit linear series}
 We investigate the Brill-Noether theory of a $2$-pointed elliptic curve (see also
\cite{EH4}),  and we prove that $\overline{\mathcal{M}}_{g,d}^r\cap \Delta _1$ is
irreducible for
$\rho(g,r,d)=-1$. 
\begin{prop}
Let $(E, p, q)$ be a two-pointed elliptic curve. Consider the
sequences
$$a:\mbox{ }a_0<a_1<\ldots a_r\leq d,\mbox{ }\mbox{ } b:\mbox{ }b_0<b_1<\ldots b_r\leq d.
$$
1. For any linear series $l=(\mathcal{L} ,V) \in G^r_d(E)$ one has that
$\rho (l, \alpha ^l(p), \alpha ^l(q))\geq -r.$ Furthermore, if $\rho (l, \alpha ^l(p), 
\alpha ^l(q))\leq -1, $ then $p-q\in \mbox{Pic}^0(E)$ is a torsion class.
\newline
2. Assume that the sequences $a$ and $b$ satisfy the inequalities:
$d-1\leq a_i+b_{r-i}\leq d,\mbox{ } i=0,\ldots ,r.$
Then there exists at most one linear series $l\in G^r_d(E)$ such that $a^l(p)
=a$ and $a^l(q)=b.$ Moreover, there exists exactly one such linear series 
$l=(\mathcal{O} _E(D), V)$ with $D\in
E^{(d)}$, if and only if for each $i=0,\ldots ,r$ the following is satisfied: if
$a_i+b_{r-i}=d$, then $D\sim a_i\ p+b_{r-i}\ q$, and if $(a_i+1)\ p+b_{r-i}\ q\sim D$, 
then $a_{i+1}=a_i+1.$
\end{prop}
\textsl{Proof:} In order to prove 1.\ it is enough to notice that for
dimensional reasons there must be sections $\sigma _i\in V$ such that
$\mbox{div}(\sigma _i)\geq a_i^l(p)\ p+a_{r-i}^l(q)\ q,$
therefore, $a_i^l(p)+a_{r-i}^l(q)\leq d.$
By adding up all these inequalities, we get that 
$\rho (l, \alpha ^l(p), \alpha ^l(q))\geq -r.$
Furthermore, $\rho (l, \alpha ^l(p),\alpha ^l(q))\leq -1$ precisely when for at
least two values $i<j$ we have equalities 
$a_i+b_{r-i}=d,\mbox{ }a_j+b_{r-j}=d,$
which means that there are sections $\sigma _i, \sigma _j\in V$ such that 
$\mbox{div}(\sigma _i)=a_i\ p+b_{r-i}\ q,\mbox{ } \mbox{div}(\sigma _j)=a_j\ p+ b_{r-j}\
q.$ By subtracting, we see that $p-q \in \mbox{Pic}^0(E)$ is torsion. 
The second part of the Proposition is in fact Prop.5.2  from \cite{EH4}.\hfill $\Box$
\vskip10pt 
\begin{prop}
Let $g,r,d$ be such that $\rho (g,r,d)=-1.$ Then the intersection $\overline
{\mathcal{M} } _{g,d}^r\cap \Delta _1$ is irreducible.
\end{prop}
\textsl{Proof: }Let $Y$ be an irreducible component of 
$\overline{\mathcal{M}}_{g,d}^r\cap \Delta _1.$ Either $Y\cap \mbox{Int}\Delta _1\neq
\emptyset,\mbox{ hence }Y=\overline{Y\cap \mbox{Int}\Delta _1} $, or $Y\subseteq \Delta
_1-\mbox{Int} \Delta _1.$  The second alternative never occurs. Indeed, if $Y\subseteq
\Delta _1-\mbox{Int}
\Delta _1,$ then since $\mbox{codim}\mbox{ } (Y,\overline{\mathcal{M} _g })=2$,  $Y$ must
be one of the irreducible components of $\Delta _1-\mbox{Int} \Delta _1.$ The components
of $\Delta _1-\mbox{Int} \Delta _1$ correspond to curves with two nodes. We list these
components (see \cite{Ed1}):
\begin{itemize}
\item For $1\leq j\leq g-2,$ $\Delta _{1j}$ is the closure of the locus in 
$\overline{\mathcal{M} } _g$ whose general point corresponds to a chain composed
of an elliptic curve, a curve of genus $g-j-1$, and a curve of genus $j$.
\item The component $\Delta _{01}$, whose general point corresponds to the
union of a smooth elliptic curve and an irreducible nodal curve of genus 
$g-2$.
\item The component $\Delta _{0,g-1}$ whose general point corresponds to the union of  a
smooth curve of genus $g-1$ and an irreducible rational curve.
\end{itemize}
As the general point of $\Delta _{1,j}, \Delta _{0,1}$ or $\Delta _{0,g-1}$ 
is a tree-like curve which satisfies the conditions of Prop.3.2 it follows that  such a
curve satisfies the `strong' Brill-Noether Theorem, hence $\Delta _{1,j}\nsubseteq
\overline{\mathcal{M}  } _{g,d}^r$, $\Delta _{0,1}\nsubseteq \overline{\mathcal{M} }
_{g,d}^r$ and
$\Delta _{0,g-1}\nsubseteq \overline{\mathcal{M} } _{g,d}^r,$ a contradiction.
So, we are left with the first possibility: $Y=\overline{Y\cap \mbox{Int} \Delta _1} .$
We are going to determine the general point $[C]\in Y\cap \mbox{Int} \Delta _1.$ 
Let $X=C\cup E, g(C)=g-1, E$ elliptic, $E\cap C=\{ p\} $ such that $X$ carries
a limit $\mathfrak g ^r_d$, say $l$. By the additivity of the Brill-Noether number, we have:
$$-1=\rho (g,r,d)\geq \rho (l,C,p)+\rho (l,E,p).$$
Since $\rho (l,E,p)\geq 0$, it follows that $\rho (l,C,p)\leq -1$, so $w^{l_C}(p)\geq r.$
Let us denote by $$\beta :\mathcal{C} _{g-1}\times \mathcal{C} _1\rightarrow \mbox{Int} \Delta _1$$
the natural map given by
$\beta ([C,p],[E,q])=[X:=C\cup E/p\sim q].$
We claim that if we choose $X$ generically, then $\alpha _0^{l_C}(p)=0.$ If not, $p$ is a base point of $l_C$ and after removing the base point we get that
$[C]\in \mathcal{M} ^r_{g-1,d-1}.$ Note that $\rho (g-1,r,d-1)=-2$, so 
$\mbox{dim } \mathcal{M}^r_{g-1,d-1}=3g-8$ (cf. \cite{Ed2}).
If we denote by $\pi :\mathcal{C} _{g-1}\rightarrow \mathcal{M} _{g-1}$ the morphism which `forgets the point', we get that
$$\mbox{dim}\mbox{ } \beta (\pi ^{-1}(\mathcal{M} ^r_{g-1,d-1})\times \mathcal{C} _1)
=3g-6< \mbox{dim}\mbox{ } Y,$$
a contradiction. Hence, for the generic $[X]\in Y$ we must have $\alpha _0^{l_C}
(p)=0$, so $a_r^{l_E}(p)=d.$ 
Since an elliptic curve cannot have a meromorphic function with a single pole,
it follows that $a_{r-1}^{l_E}(p)\leq d-2$ and this implies 
$ \alpha ^{l_C}(p)\geq (0,1,\ldots ,1),$
i.e. $l_C$ has a cusp at $p$. Thus, if we introduce the notation 
$$\mathcal{C} ^r_{g-1,d}(0,1,\ldots ,1)=\{ [C,p]\in \mathcal{C} _{g-1}:G^r_d
(C, (p,(0,1,\ldots ,1)))\neq \emptyset \} ,$$
then $Y\subseteq \overline{\beta (\mathcal{C} ^r_{g-1,d}(0,1.\ldots ,1)\times
\mathcal{C} _1)}.$
On the other hand, it is known (cf. \cite{EH2}) that $\mathcal{C} _{g-1,d}^r(0,1,\ldots ,1)$ is irreducible of dimension $3g-6$ (that is, codimension $1$ in $\mathcal{C} _{g-1}$), so we must have
$Y=\overline{\beta (\mathcal{C}_{g-1,d}^r(0,1,\ldots ,1)\times \mathcal{C} _1)}$,  which
not only proves that $\overline{\mathcal{M} }_{g,d}^r\cap \Delta _1$ is irreducible, but
also determines the intersection.\hfill $\Box$
\section{The Kodaira dimension of $\mathcal{M}_{23}$ }
In this section we prove that $\kappa (\mathcal{M}_{23})\geq 2$ and we investigate closely
the multicanonical linear systems on $\overline{\mathcal{M} }_{23}$. We now describe the
three multicanonical Brill-Noether divisors from Section 2.
\subsection{The divisor $\overline{\mathcal{M} }_{12}^1$}
There is a stratification of $\mathcal{M}_{23}$ given by gonality:
$$\mathcal{M}^1_2\subseteq \mathcal{M}^1_3\subseteq  \ldots \subseteq 
\mathcal{M}_{12}^1\subseteq \mathcal{M}_{23}.$$
For $2\leq d\leq g/2+1$ one knows that $\mathcal{M}^1_k=\mathcal{M}_{g,k}^1$ is an
irreducible variety of dimension $2g+2d-5$. The general point of $\mathcal{M}_{g,d}^1$
corresponds to a curve having a unique $\mathfrak g ^1_d.$
%%
%%%
%%
\subsection{The divisor $\overline{\mathcal{M}}_{17}^2$}
The {\sl Severi variety } $V_{d,g}$ of irreducible plane curves of degree $d$
and geometric genus $g$, where $0\leq g\leq {d-1\choose 2}$, is an irreducible
subscheme of ${\mathbb P}^{d(d+3)/2}$ of dimension $3d+g-1$ (cf. \cite{H},
\cite{Mod}). Inside $V_{d,g}$ we consider the open dense subset $U_{d,g} $ of
irreducible plane curves of degree $d$ having exactly $\delta ={d-1\choose 2} -g$
nodes and no other singularities. There is a global normalization map 
$$m:U_{d,g}\rightarrow \mathcal{M} _g,\mbox{ } m([Y]):=[\tilde{Y}],\mbox{ }
 \tilde{Y}\mbox{ is the normalization of } Y. $$
When $d-2\leq g\leq {d-1\choose 2}, d\geq 5,$ $U_{d,g}$ has the expected
number of moduli, i.e.
$$\mbox{dim}\mbox{ } m(U_{d,g})=\mbox{min}(3g-3,3g-3+\rho(g,2,d)).$$
In our case we can summarize this as follows:
\begin{prop}
There is exactly one component of $\mathcal{G}^2_{17}$  mapping dominantly to
$\mathcal{M} _{17}^2$. The general element $(C,l)\in \mathcal{G}^2_ {17}$ corresponds
to a curve $C$ of genus $23$,  together with a $\mathfrak g^2_{17}$ which provides a
plane model for $C$ of degree $17$ with $97$ nodes.
\end{prop}
\subsection{The divisor $\overline{\mathcal{M} } _{20}^3$}
Here we combine the result of Eisenbud and Harris (see \cite{EH2})
about the uniqueness of divisorial components of $\mathcal{G} ^r_d$
when $\rho (g,r,d)=-1$, with Sernesi's (see \cite{Se2}) which asserts the existence of
components of the Hilbert scheme $H_{d,g}$ parametrizing curves in $\mathbb P^3$ of degree
$d$ and genus $g$ with the expected number of moduli, for $d-3\leq g\leq 3d-18,d\geq 9$.
\begin{prop}
There is exactly one component of $\mathcal{G} ^3_{20}$ mapping dominantly to 
$\mathcal{M} _{20}^3$.  The general point of this component corresponds to a pair $(C, l)$
where $C$ is a curve of genus $23$ and $l$ is a very ample $\mathfrak g ^{3}_{20}$.
\end{prop}
We are going to prove that the Brill-Noether divisors $\overline{\mathcal{M} }
^1_{12},\overline{\mathcal{M} } _{17}^2$ and $\overline{\mathcal{M} }^3_{20}$
are mutually distinct.
\begin{thm}
There exists a smooth curve of genus $23$ having a $\mathfrak g^2_{17},$ but no
$\mathfrak g^3_{20}$'s. Equivalently, one has ${\rm supp} (\mathcal{M} ^2_{17}) 
\nsubseteq {\rm supp} (\mathcal{M} ^3_{20}).$
\end{thm}
{\sl Proof: } It suffices to construct a reducible curve $X$ of compact type of genus $23$, which has a smoothable limit $\mathfrak g^2_{17}$, but no limit $\mathfrak g^3_{20}.$ If $[C]\in  \mathcal{M} _{23}$ is a nearby smoothing of $X$ which preserves the $\mathfrak g^2_{17}$, then $[C]\in \mathcal{M} ^2_{17}-
\mathcal{M} ^3_{20}.$
Let us consider the following curve:\newline
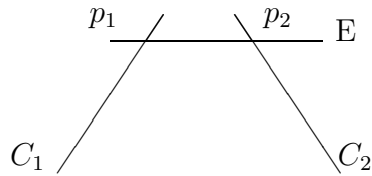
\begin{picture}(300,80)
\put(165,10){\line(2,3){40}}
\put(147,13){$C_1$}
\put(185,60){\line(1,0){80}}
\put(177,67){$p_1$}
\put(243,67){$p_2$}
\put(272,10){\line(-2,3){40}}
\put(271,13){$C_2$}
\put(270,60){E}
\end{picture}

$$X:=C_1\cup C_2\cup E,$$
where $(C_1,p_1)$ and $(C_2,p_2)$ are general pointed curves of genus $11$,
$E$ is an elliptic curve, and $p_1-p_2$ is a primitive $9$-torsion point in 
$\mbox{Pic}^0(E)$\newline {\textbf {Step 1)} } \textit{There is no limit $\mathfrak
g^3_{20}$ on $X$.} Assume that $l$ is a limit $\mathfrak g^3_{20}$ on $X$. By 
the additivity of the Brill-Noether number, 
$$-1\geq \rho (l_{C_1},p_1)+\rho (l_{C_2},p_2)+\rho (l_E, p_1, p_2).$$
Since $(C_i,p_i)$ are general points in $\mathcal{C} _{11}$, it follows from 
Prop.3.2 that $\rho (l_{C_i},p_i)\geq 0,$ hence $\rho (l_E,p_1,p_2)\leq -1.$
On the other hand $\rho (l_E,p_1,p_2)\geq -3$ from Prop.4.1.
\newline
\indent Denote by $(a_0,a_1,a_2,a_3)$ the vanishing sequence of $l_E$ at $p_1$, and by 
$(b_0,b_1,b_2,b_3)$ that of $l_E$ at $p_2.$ 
The condition (8) for a general pointed curve $[(C_i,p_i)]\in \mathcal{C}_{11}$ to 
possess a $\mathfrak g^3_{20}$ with prescribed ramification at the point $p_i$ and the
compatibility conditions between $l_{C_i}$ and $l_E$ at
$p_i$ give that:
\begin{equation}
(14-a_3)_+ +(13-a_2)_+ +(12-a_1)_+ +(11-a_0)_+\leq 11,
\end{equation}
and
\begin{equation}
(14-b_3)_{+} +(13-b_2)_{+} +(12-b_1)_+ +(11-b_0)_+\leq 11.
\end{equation}
{\sl 1st case: } $\rho(l_E,p_1,p_2)=-3.$ Then $a_i+b_{3-i}=20$, for $i=0,\ldots ,3$
and it immediately follows that $20(p_1-p_2)\sim 0$ in $\mbox{Pic}^0(E)$, a 
contradiction.
\newline
{\sl 2nd case: } $\rho(l_E,p_1,p_2)=-2.$ We have two distinct possibilities here: 
i) $a_0+b_3=20,a_1+b_2=20,a_2+b_1=20,a_3+b_0=19.$ Then it follows that
$a^{l_E}(p_1)=(0,9,18,19)$ and $a^{l_E}(p_2)=(0,2,11,20)$, while according to (9),
$a_3\leq 15$, (because $\rho (l_{C_1},p_1)\leq 1)$, a contradiction. ii)
$a_0+b_3=20,a_1+b_2=20,a_2+b_1=19,a_3+b_0=20$. Again, it follows that $a_3=a_0+18\geq 15$,
a contradiction.
\newline
{\sl 3rd case: } $\rho(l_E,p_1,p_2)=-1.$ Then $\rho(l_{C_i},p_i)=0$ and $l$ is
a refined limit $\mathfrak g^3_{20}$.  From (9) and (10) we must have:
$a^{l_E}(p_i)\leq (11,12,13,14),\mbox{ } i=1,2.$ There are four possibilities:
i) $a_0+b_3=a_1+b_2=20, a_2+b_1=a_3+b_0=19.$ Then $a_1=a_0+9\leq 12,$
so $b_3=20-a_0\geq 17,$ a contradiction. ii) $a_0+b_3=a_2+b_1=20,a_2+b_1=a_3+b_0=19.$ 
Then $b_3=20-a_0\leq 14,$ so $a_2=a_0+9\geq 15,$ a contradiction.
iii) $a_0+b_3=a_3+b_0=20,a_1+b_2=a_2+b_1=19.$ Then $b_3=19-a_0\leq 14,$
so $a_3\geq a_0+9\geq 15$, a contradiction.
iv) $a_0+b_3=a_3+b_0=19,a_1+b_2=a_2+b_1=20.$ Then $b_3=19-a_0\leq 14,$
so $a_2\geq a_1+9\geq 15$, a contradiction again.
We conclude that $X$ has no limit $\mathfrak g^3_{20}.$
\newline
\newline 
\textbf{Step 2) } \textit{There exists a smoothable limit $\mathfrak g^2_{17}$ on $X$,
hence $[X]\in \overline{\mathcal{M} } ^2_{17}.$} We construct a limit linear  series $l$
of type $\mathfrak g^2_{17}$ on $X$, aspect by aspect: on $C_i$ take $l_{C_i}\in
G^2_{17}(C_i)$ such that $a^{l_{C_i}}(p_i)= (4,9,13).$ Note that in this case
$\sum_{j=0}^r(\alpha _j+g-d+r)_+=g,$ so (8) ensures the existence of such a  $\mathfrak
g^2_{17}.$ On $E$ we take $l_E=|V_E|$, where $|V_E|\subseteq |4p_1+13p_2|=|4p_2+13p_1|$
is a $\mathfrak g^2_{17}$ with vanishing sequence $(4,8,13)$ at $p_i.$ Prop.4.1 ensures the existence of such a linear series.
In this way $l$ is a refined limit $\mathfrak g^2_{17}$ on $X$ with $\rho(
l_{C_i},p_i)=0,\rho(l_E,p_1,p_2)=-1.$ We prove that $l$ is dimensionally proper.
Let 
$\pi _i:\mathcal{C} _i\rightarrow \Delta _i,\mbox{ } \tilde p_i:\Delta _i\rightarrow 
\mathcal{C} _i,$ be the versal deformation of $[(C_i,p_i)]\in \mathcal{C} _{11},$ and
$\sigma _i:
\mathcal{G}^2_{17}(\mathcal{C}_i/\Delta _i, (\tilde{p}_i,(4,8,11)))\rightarrow \Delta_i$
the projection.
\newline 
\indent Since being general is an open condition, we have that $\sigma _i$ is surjective 
and $\mbox{dim}\mbox{ } \sigma_i ^{-1}(t)=\rho (l_{C_i},p_i)=0$, for each $t\in \Delta
_i$, therefore
$$\mbox{dim}\mbox{ } \mathcal{G} ^2_{17}(\mathcal{C} _i/\Delta _i, (\tilde
p_i,(4,8,11)))=\mbox{dim}\mbox{ } \Delta _i+\rho (l_{C_i},p_i)=31.$$
\newline
Next, let
$\pi:\mathcal{C}\rightarrow \Delta,\mbox{ } \tilde p_1,\tilde p_2: \Delta \rightarrow
\mathcal{C} $ be the versal deformation of $(E,p_1,p_2).$ We prove that 
$$\mbox{dim } \mathcal{G} ^2_{17}(\mathcal{C}/\Delta, (\tilde p_i,(4,7,11)))=
\mbox{dim } \Delta +\rho (l_E,p_1,p_2)=1.$$
This follows from Prop.4.1, since a $2$-pointed elliptic curve $(E_t,\tilde 
p_1(t),\tilde p_2(t))$ has at most one $\mathfrak g^2_{17}$ with ramification
$(4,7,11)$ at both $\tilde p_1(t)$ and $\tilde p_2(t)$, and exactly one when
$9(\tilde p_1(t)-\tilde p_2(t))\sim 0.$ 
Hence $\mbox{Im} \mathcal{G}^2_{17}(\mathcal{C}/\Delta, (\tilde{p}_i,(4,7,11))) =\{t\in
\Delta:9(\tilde p_1(t)-\tilde p_2(t))\sim 0\mbox{ in } \mbox{Pic}^0(E_t)\}$, which is a
divisor on $\Delta $, so the claim follows and $l$ is a dimensionally proper $\mathfrak
g^2_{17}.$\hfill $\Box $
\vskip10pt
A slight variation of the previous argument gives us:
\begin{prop}
We have ${\rm supp} (\overline{\mathcal{M} } ^2_{17}\cap \Delta _1)\neq {\rm supp } 
(\overline{\mathcal{M} } ^3_{20}\cap \Delta _1).$
\end{prop}
{\textsl Proof:} We construct a curve $[Y]\in 
\Delta _1\subseteq \overline{\mathcal{M}}_{23}$ which has a smoothable limit $\mathfrak
g^2_{17}$ but no limit $\mathfrak g^3_{20}.$ Let us consider the following curve:\newline
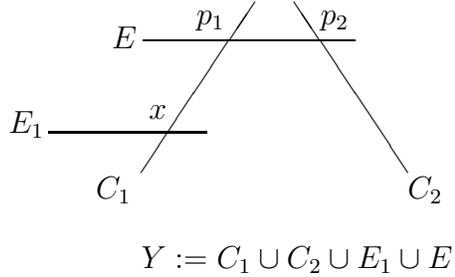
\begin{picture}(300,80)
\put(132,25){\line(1,0){60}}
\put(167,10){\line(2,3){43}}
\put(168,60){\line(1,0){80}}
\put(268,10){\line(-2,3){43}}
\put(170,30){$x$}
\put(188,66){$p_1$}
\put(235,66){$p_2$}
\put(117,25){$E_1$}
\put(150,0){$C_1$}
\put(268,0){$C_2$}
\put(156,57){$E$}
\end{picture}
%%
%%
%%%
$$Y:=C_1\cup C_2\cup E_1\cup E,$$
where $(C_2,p_2)$ is a general point of $\mathcal{C} _{11}$, $(C_1,p_1,x)$ is
a general $2$-pointed curve of genus $10$, $(E_1,x)$ is general in $\mathcal{C}_1$, $E$ is
an elliptic curve, and $p_1-p_2\in \mbox{Pic}^0(E)$ is a primitive $9$-torsion. In order
to prove that $Y$ has no limit $\mathfrak g^3_{20}$, one just has to take into account
that according to Prop.3.3, the condition for a general $1$-pointed curve $(C,z)$ of genus
$g$, to have a $\mathfrak g^r_d$ with ramification $\alpha $ at $z$ is the same with the
condition for a general $2$-pointed curve $(D,x,y)$ of genus $g-1$ to have a $\mathfrak
g^r_d$ with ramification $\alpha $ at $x$ and a cusp at $y$. Therefore we can repeat what
we did in the proof of Theorem 2. Next, we construct $l$, a smoothable limit $\mathfrak
g^2_{17}$ on $Y$: take $l_{C_2}\in G^2_{17}(C_2,(p_2,(4,8,11))),l_E=|V_E|\subseteq
|4p_1+13p_2|$, with $\alpha ^{l_E}(p_i)=(4,7,11),$ on $E_1$ take $l_{E_1}=14x+|3x|$, and
finally on $C_1$ take $l_{C_1}$ such that $\alpha ^{l_{C_1}}(p_1)=(4,8,11),\alpha
^{l_{C_1}}(x)=(0,0,1).$ Prop.3.3 ensures the existence of $l_{C_1}.$ Clearly, $l$ is a
refined limit $\mathfrak g^2_{17}$ and the proof that it is  smoothable is all but
identical to the one in the last part of Theorem 2.\hfill
$\Box $
\newline
\indent The other cases are settled by the following:
\begin{thm}
There exists a smooth curve of genus $23$ having a $\mathfrak g^1_{12}$ but  having no
$\mathfrak g^2_{17}$ nor $\mathfrak g^3_{20}$. Equivalently, ${\rm supp} (\mathcal{M}
^1_{12})\nsubseteq  {\rm supp} (\mathcal{M} ^2_{17})$ and
${\rm supp}(\mathcal{M} ^1_{12})\nsubseteq {\rm supp}(\mathcal{M}^3_{20}).$
\end{thm}
%%\eject
\textsl{Proof: } We take the curve considered in \cite{EH3}:
\newline
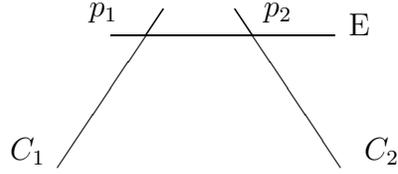
\begin{picture}(300,80)
\put(165,10){\line(2,3){40}}
\put(147,13){$C_1$}
\put(185,60){\line(1,0){85}}
\put(177,67){$p_1$}
\put(243,67){$p_2$}
\put(272,10){\line(-2,3){40}}
\put(281,13){$C_2$}
\put(275,60){E}
\end{picture}

$$Y:=C_1\cup C_2\cup E,$$
where $(C_i,p_i)$ are general points of $\mathcal{C} _{11}$, $E$ is elliptic
and $p_1-p_2\in \mbox{Pic} ^0(E)$ is a primitive $12$-torsion.
Clearly $Y$ has a (smoothable) limit $\mathfrak g^1_{12}$: on $C_i$ take the
pencil $|12p_i|$, while on $E$ take the pencil spanned by $12p_1$ and $12p_2.$
It is proved in \cite{EH3} that $Y$ has no limit $\mathfrak g^2_{17}$'s and similarly
one can prove that $Y$ has no limit $\mathfrak g^3_{20}$'s either. We omit the
details.\hfill $\Box $
\newline
\newline
\indent Now we are going to prove that equation (2)
$$\mbox{supp}(\mathcal{M} ^1_{12})\cap \mbox{supp}(\mathcal{M} ^2 _{17})=
\mbox{supp}(\mathcal{M}  ^2_{17})\cap \mbox{supp}(\mathcal{M} ^3_{20})=
\mbox{supp}(\mathcal{M} ^3_{20})\cap \mbox{supp}(\mathcal{M} ^1_{12})
$$
 is impossible, and as explained before, this will imply that $\kappa
(\mathcal{M}_{23})\geq 2.$ The main step in this direction is the following:
\begin{prop}
There exists a stable curve of compact type of genus $23$ which has  a smoothable limit
$\mathfrak g^3_{20}$, a smoothable limit $\mathfrak g^2_{15}$ (therefore also a $\mathfrak
g^2_{17}$), but has generic gonality, that is, it does not have any limit $\mathfrak
g^1_{12}$.
\end{prop}
%%\eject
\textsl{ Proof }  We shall consider the following stable curve $X$ of genus
 $23$:\newline
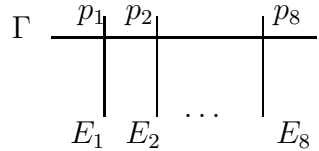
\begin{picture}(300,60)
\put(160,40){\line(1,0){100}}
\put(180,10){\line(0,1){38}}
\put(200,10){\line(0,1){38}}
\put(240,10){\line(0,1){38}}
\put(210,10){$\ldots $}
\put(145,40){$\Gamma $}
\put(167,0){$E_1$}
\put(188,0){$E_2$}
\put(245,0){$E_8$}
\put(170,47){$p_1$}
\put(188,47){$p_2$}
\put(244,47){$p_8$}
\end{picture}
$$
X:=\Gamma \cup E_1\cup \ldots \cup E_8,
$$
where the $E_i$'s are elliptic curves, $\Gamma \subseteq \mathbb P^2$ is a 
general smooth plane septic and the points of attachment $\{p_i\} =\Gamma \cup
E_i$ are general points of $\Gamma $.
\newline
\textbf{Step 1)} $\textit{\mbox{There is no limit }$\mathfrak g^1_{12}$\mbox{ on  }X}$. 
Assume that $l$
is a limit $\mathfrak g^1_{12}$ on $X$. Since the ellip\-tic curves $E_i$ cannot have
 meromorphic functions with a single pole, we have that $a^{l_{E_i}}(p_i)\leq (10,12)$,
hence $\alpha ^{l_{\Gamma}}(p_i)\geq (0,1),$ that is, $l_{\Gamma }$ has a cusp at $p_i$
for $i=1,\ldots ,8.$ We now prove that $\Gamma $ has no $\mathfrak g^1_{12}$'s with cusps
at the points $p_i.$
\newline
\indent First, we notice that $\mbox{dim }G^1_{12}(\Gamma )=\rho(15,1,12)=7.$ Indeed, 
if we assume that $\mbox{dim }G^1_{12}(\Gamma) \geq 8$, by applying Keem's Theorem
(cf.\
\cite{ACGH}, p.200) we would get that $\Gamma $ possesses a $\mathfrak g^1_4$,  which
is impossible since $\mbox{gon}(\Gamma )=6.$ (In general, if $Y\subseteq \mathbb P^2$
is a smooth plane curve, $\mbox{deg}(Y)=d,$ then $\mbox{gon }(Y)=d-1$, and the
$\mathfrak g^1_{d-1}$ computing the gonality is cut out by the lines passing through a
point 
$p\in Y,$ see \cite{ACGH}.) Next, we define the variety
$$\Sigma =\{(l,q_1,\ldots ,q_8)\in G^1_{12}(\Gamma )\times  \Gamma ^8:\alpha
^l(q_i)\geq  (0,1),i=1,\ldots ,8\}$$
and denote by $\pi_1:\Sigma \rightarrow G^1_{12}(\Gamma )$ and $\pi_2: \Sigma
\rightarrow 
\Gamma ^8$ the two projections. For any $l\in G^1_{12}(\Gamma )$,
the fibre $\pi_1^{-1}(l)$ is finite, hence $\mbox{dim }\Sigma =\mbox{dim
}G^1_{12}(\Gamma)=7$, which shows that $\pi_2$ cannot be surjective and this proves our
claim.
\newline
\textbf{Step 2)} \textit{\mbox{There exists a smoothable limit }$\mathfrak g^2_{17}$
\mbox { on }$X$,\mbox{ hence }$[X]\in \overline{\mathcal{M} } _{17}^2$.}
We construct $l$, a limit $\mathfrak g^2_{17}$ on $X$ as follows: 
on $\Gamma $ there is a (unique) $\mathfrak g^2_7$, and we consider $l_{\Gamma
}=\mathfrak g^2_7(p_1+\cdots +p_8)$, i.e. the $\Gamma -$ aspect $l_{\Gamma }$ is
obtained from the $\mathfrak g^2_7$ by adding the base points $p_1,\ldots ,p_8$.
Clearly $a^{l_{\Gamma }}(p_i)=(1, 2, 3)$ for each $i$. On $E_i$ we take
$l_{E_i}=\mathfrak g^2_3(12p_i)$ for $i=1,\ldots ,8,$ where $\mathfrak g^2_3$ is a
complete linear series of the form $|2p_i+x_i|$, with $x_i\in E_i-\{ p_i\}.$
Furthermore, $a^{l_{E_i}}(p_i)=(12, 13, 14)$, so $l=\{l_{\Gamma },l_{E_i}\}$ is a
refined limit $\mathfrak g^2_{15}$ on $X$. One sees that $\rho (l_{E_i},
\alpha ^{l_{E_i}}(p_i))=1$ for all $i$, 
$\rho(l_{\Gamma },\alpha ^{l_{\Gamma }}(p_1),\ldots ,\alpha ^{l_{\Gamma }}(p_8))=-15,$
and $\rho(l)=-7.$ We now prove that $l$ is dimensionally proper.\newline
\indent Let 
$\pi _i:\mathcal{C} _i\rightarrow \Delta _i,\mbox{ } \tilde{p} _i:\Delta _i\rightarrow
\mathcal{C} _i$ be the versal deformation space of $(E_i,p_i),$ for $i=1,\ldots ,8.$
There is an obvious isomorphism over $\Delta _i$
$$\mathcal{G} ^2_{15}(\mathcal{C} _i/\Delta _i,(\tilde{p} _i,(12,12,12)))\simeq 
\mathcal{G} ^2_3(\mathcal{C} _i/\Delta _i,(\tilde{p} _i,0)).$$ If $\sigma
_i:\mathcal{G} ^2_3(\mathcal{C} _i/\Delta _i,(\tilde{p} _i,0))\rightarrow \Delta _i$
is the natural projection, then for each $t\in \Delta _i$, the fibre $\sigma
_i^{-1}(t)$ is isomorphic to $\pi _i ^{-1}(t)$, the isomorphism being given by 
$$\pi _i^{-1}(t)\ni q\mapsto |2\tilde{p} _i(t)+q|\in G^2_3(\pi _i^{-1}(t)).$$
Thus, $\mathcal{G} ^2_3(\mathcal{C} _i/\Delta _i,(\tilde{p}_i,0))$ is a smooth 
irreducible surface, which shows that $l$ is dimensionally proper w.r.t.$\ E_i$. Next,
let us consider $\pi :\mathcal{X} \rightarrow \Delta,\mbox{ } \tilde{p}_1,\ldots
,\tilde{p} _8:\Delta \rightarrow \mathcal{X}$, the versal deformation of $(\Gamma,
p_1,\ldots ,p_8).$ We have to prove that 
$$\mbox{dim }\mathcal{G}^2_{15}(\mathcal{X} /\Delta,(\tilde{p} _i,(1,1,1)))=
\mbox{dim }\Delta +\rho(l_{\Gamma },\alpha ^{l_{\Gamma}}(p_i))=35.$$
There is an isomorphism over $\Delta $,
$$\mathcal{G}^2_{15}(\mathcal{X}/\Delta ,(\tilde{p}_i,(1,1,1)))\simeq
\mathcal{G}^2_{7}(\mathcal{X}/\Delta, (\tilde{p}_i,0)).$$
If $\pi_0:\mathcal{C}\rightarrow \mathcal{M}$ is the versal deformation space of
$\Gamma $, then we denote by $\mathcal{G}^2_7\rightarrow \mathcal{M} $ the scheme 
parametrizing $\mathfrak g^2_7$'s on curves of genus $15$ `nearby' $\Gamma$ (See
Section 3 for this notation). Clearly
$\mathcal{G}^2_7(\mathcal{X}/\Delta,(\tilde{p}_i,0))\simeq \mathcal{G}^2_7\times
_{\mathcal{M} }\Delta ,$ so it suffices to prove that $\mathcal{G}^2_7$ has the
expected dimension at the point $(\Gamma, \mathfrak g^2_7).$ For this we use Prop.
3.1. We have that 
$N_{\Gamma /\mathbb P^2}=\mathcal{O}_{\Gamma }(7), K_{\Gamma }=\mathcal{O}_{\Gamma
}(4)$, hence $$H^1(\Gamma, N_{\Gamma /\mathbb P^2})\simeq H^0(\Gamma
,\mathcal{O}_{\Gamma }(-3))^{\vee}=0,$$ so $l$ is dimensionally proper w.r.t.\ 
$\Gamma$ as well. We conclude that $l$ is smoothable.\newline
\newline
\textbf{Step 3)} \textit{There exists a smoothable limit $\mathfrak g^3_{20}$ on  $X$,
that is $[X]\in \overline{\mathcal{M}}  ^3_{20}$.} First we notice that there is an
isomorphism $\Gamma \xrightarrow{\sim } G^1_6(\Gamma )$, given by $$\Gamma \ni
p\mapsto |\mathfrak g^2_7-p|\in G^1_6(\Gamma ).$$ Consequently, there is a
$2$-dimensional family of $\mathfrak g^3_{12}$'s on
$\Gamma $, of the form  $\mathfrak g^3_{12}=\mathfrak g^1_6+\mathfrak
h^1_6=|2\mathfrak g^2_7-p-q|$, where $p,q\in \Gamma $. Pick $l_0=l_0'+l_0''$, with
$l_0',l_0''\in G^1_6(\Gamma )$, a general $\mathfrak g^3_{12}$ of this type. We
construct $l$, a limit $\mathfrak g^3_{20}$ on $X$, as follows: the $\Gamma $-aspect
is given by $l_{\Gamma }=l_0(p_1+\cdots p_8),$ and because of the generality of the
chosen $l_0$ we have that $\rho(l_{\Gamma },\alpha ^{l_{\Gamma }}(p_1),\ldots ,\alpha
^{l_{\Gamma }}(p_8))=-9.$ The $E_i$-aspect is given by $l_{E_i}=\mathfrak
g^3_4(16p_i),$ where $\mathfrak g^3_4=|3p_i+x_i|$, with $x_i\in E_i-\{p_i\}$, for
$i=1,\ldots ,8.$ It is clear that $\rho (l_{E_i},\alpha ^{l_{E_i}}(p_i))=1$ and that
$l'=\{l_{\Gamma },l_{E_i}\}$ is a refined limit $\mathfrak g^3_{20}$ on $X$.
\newline
\indent In order to prove that $l'$ is dimensionally proper, we first notice 
that $l'$ is dimensionally proper w.r.t.$\ $the elliptic tails $E_i$. We now prove that
$l'$ is dimensionally proper w.r.t.$\ \Gamma $. As in the previous step, we consider
$\pi:\mathcal{X}\rightarrow \Delta ,\mbox{ } \tilde{p} _1,\ldots ,\tilde{p} _8:\Delta
\rightarrow \mathcal{X},$ the versal deformation of $(\Gamma ,p_1,\ldots ,p_8)$ and
$\pi_0:\mathcal{C}\rightarrow \mathcal{M}$, the versal deformation space of $\Gamma $.
There is an isomorphism over $\Delta $
$$\mathcal{G} ^3_{20}(\mathcal{X}/\Delta , (\tilde{p} _1,\alpha ^{l_{\Gamma
}}(p_1),\ldots ,(\tilde{p} _8,\alpha ^{l_{\Gamma }}(p_8)))\simeq
\mathcal{G}^3_{12}(\mathcal{C}/\mathcal{M})\times _{\mathcal{M} }\Delta.$$ It suffices
to prove that $\mathcal{G} ^3_{12}=\mathcal{G} ^3_{12}(\mathcal{C}/\mathcal{M})$ has a
component of the expected dimension passing through $(\Gamma ,l_0)$. In this way, the
genus $23$ problem is turned into a deformation theoretic problem in genus $15$.
Denote as usual by $\sigma :
\mathcal{G}^3_{12}\rightarrow \mathcal{M}$ the natural projection.  According to
Prop.3.1, it will be enough to exhibit an element $(C,l)\in \mathcal{G} ^3_{20}$,
sitting in the same component as $(\Gamma ,l_0)$, such that the linear system $l$ is
base point free and simple, and if $\phi :C\rightarrow \mathbb P^3$ is the map induced
by $l$, then $H^1(C,N_{\phi })=0$. Certainly we cannot take $C$ to be a smooth plane
septic because in this case $H^1(C, N_{\phi })\neq 0$, as one can easily see. Instead,
we consider the $6$-gonal locus in a neighbourhood of the point $[\Gamma ]\in
\mathcal{M} _{15}$, or equivalently, the $6$-gonal locus in $\mathcal{M} $, the versal
deformation space of $\Gamma $. One has the projection $\mathcal{G} ^1_6\rightarrow
\mathcal{M}$, and the scheme $\mathcal{G} ^1_6$ is smooth (and irreducible) of
dimension $37(=2g+2d-5; g=15,d=6)$. We denote by 
$$\mu :\mathcal{G} ^1_6\times _{\mathcal{M} }\mathcal{G}^1_6\rightarrow \mathcal{M}
,\mbox{ }\mbox{ } \mu ([C, l,l'])=[C] .$$ Let us pick a component $\mathcal{X}
\subseteq \mathcal{G}^1_6\times _{\mathcal{M} }\mathcal{G}^1_6$ such that $(\Gamma
,l_0', l_0'')\in \mathcal{X}.$ The general point  of $\mathcal{X}$ corresponds to a
curve $C$ with two base-point-free pencils $l',l''\in G^1_6(C)$ such that if
$f':C\rightarrow \mathbb P^1$ and $f'':C\rightarrow \mathbb P^1$ are the corresponding
morphisms, then 
$$\phi =(f',f''):C\rightarrow \mathbb P^1\times \mathbb P^1$$
is birational. We denote by $\eta :\mathcal{X}\rightarrow \mathcal{G}^3_{12}$, the map
given by $\eta (C, l', l'') :=(C, l'+l'')$. The fact that $\eta $ is well-defined
follows from the base-point-free-pencil-trick. There is a stratification of
$\mathcal{M} $ given by the number of pencils: for $i\geq 0$ we define,
$$\mathcal{M} (i)^0:=\{ [C]\in \mathcal{M} :C \hbox{ possesses } i  \hbox{ mutually
independent, base-point-free }  \mathfrak g _6^1\mbox{'s }  \},$$ and  $\mathcal{M}
(i): =\overline{\mathcal{M} (i)^0} $. The strata  $\mathcal{M} (i)^0$ are 
constructible subsets of $\mathcal{M} $, the first stratum $\mathcal{M} (1)=\mbox{Im
}(\mathcal{G} ^1_6)$ is just the $6$-gonal locus; the stratum $\mathcal{M}(2)$ is
irreducible and
$\mbox{dim }\mathcal{M}(2)=g+4d-7=32$ (cf. \cite {AC1}). We denote by $\mathcal{M}
_{sept}:=\overline{m(U_{7,15})\cap \mathcal{M}}$, the closure of the locus of smooth
plane septics in $\mathcal{M}$, and by $\mathcal{M}_{oct}:=\overline{m(U_{8,15})\cap
\mathcal{M}}$, the locus of curves which are normalizations of plane octics with $6$
nodes. Since the Severi varieties $U_{7,15}$ and $U_{8,15}$ are irreducible, so are
the loci $\mathcal{M}_{sept}$ and $\mathcal{M}_{oct}$. Furthermore $\mbox{dim
}\mathcal{M} _{sept}=27$ and $\mbox{dim }\mathcal{M} _{oct}=30$. We prove that
$\mathcal{M} _{sept}\subseteq \mathcal{M}_{oct}$, hence $\mathcal{M}_{oct}\subseteq
\mu (\mathcal{X}).$ Indeed, let us pick $Y\subseteq \mathbb P^2$ a smooth plane
septic, and $L\subseteq \mathbb P^2$ a general line, $L.Y=p_1+\cdots +p_7$. Denote
$X:=C\cup L$, $\mbox{deg }(X)=8,p_a(X)=21$. We consider the node $p_7$ unassigned,
while 
$p_1, \ldots p_6$ are assigned. By using \cite{Ta} Theorem 2.13, there exists
a flat family of plane curves $\pi :\mathcal{X}\rightarrow B$ and a point $0\in B$, 
such that $X_0=\pi ^{-1}(0)=X$, while for $0\neq b\in B$, the fibre $X_b$ is an
irreducible octic with nodes $p_1(b), \ldots p_6(b)$, and such that $p_i(b)\rightarrow
p_i$, when $b\rightarrow 0$, for $i=1,\ldots ,6$. If $\mathcal{X}'\rightarrow B$ is
the family resulting by normalizing the surface $\mathcal{X}$, and $\eta
:\mathcal{X}''\rightarrow B$ is the stable family associated to the semistable family
$\mathcal{X}'\rightarrow B$, then we get that $\eta ^{-1}(0)=Y$, while $\eta^{-1}(b)$
is the normalization of $X_b$ for $b\neq 0$. This proves our contention.  
\newline
\indent We are going to show that given a general point $[C]\in \mathcal{M}_{oct}$, 
and
$(C,l,l')\in \mu ^{-1}([C])$, the condition $H^1(C, N_{\phi })=0$ is satisfied,
hence  $\mathcal{G} ^3_{12}$ is smooth of the expected dimension at the point 
$(C, l+l')$. This will prove the existence of a component of $\mathcal{G} ^3_{12}$ 
passing through $(\Gamma, l_0)$ and having the expected dimension. We take
$\overline{C}\subseteq \mathbb P^2$, a general point of $U_{8,15}$, with nodes
$p_1,\ldots ,p_6\in \mathbb P^2$ in general position. Theorem 3.2 from $\cite{AC1}$
ensures that there exists a plane octic having 6 prescribed nodes in general position.
Let $\nu:C\rightarrow \overline{C}$ be the normalization map, $\nu
^{-1}(p_i)=q_i'+q_i''$ for $i=1,\ldots ,6$. Choose two nodes, say $p_1$ and $p_2$, and
denote by $\mathfrak g^1_6=|H-q_1'-q_1''|$ and $\mathfrak h^1_6=|H-q_2'-q_2''|$, the
linear series obtained by projecting $\overline{C} $ from $p_1$ and $p_2$
respectively. Here $H\in |\nu ^*\mathcal{O} _{\mathbb P^2}(1)|$ is an arbitrary line
section of $C$. The morphism induced by $(\mathfrak g^1_6, \mathfrak h^1_6)$ is
denoted by $\phi :C\rightarrow \mathbb P^1\times \mathbb P^1$, and  $\phi _1=s\circ
\phi:C\rightarrow \mathbb P^3$, with $s:\mathbb P^1\times \mathbb P^1\rightarrow
\mathbb P^3$ the Segre embedding. There is an exact sequence over $C$ 
\begin{equation}
0\longrightarrow N_{\phi }\longrightarrow N_{\phi _1}\longrightarrow \phi^*N_{\mathbb
P^1\times \mathbb P^1/\mathbb P^3} \longrightarrow 0.
\end{equation}
We can argue as in \cite{AC2} p.473, that for a general $(C,\mathfrak g^1_6, \mathfrak
h^1_6)$ with $[C]\in \mathcal{M}_{oct}$, we have $h^1(C, N_{\phi })=0$. Indeed, let us
denote by $\mathcal{X}_0$ the open set of $\mathcal{X}$ corresponding to points
$(X,l,l')$ such that $\chi :X\rightarrow \mathbb P^1\times \mathbb P^1$, the morphism
associated to the pair of pencils $(l,l')$ is birational, and by $\mathcal{U}
\subseteq \mathcal{X}_0$ the variety of those points $(X,l,l')\in \mathcal{X}_0 $ such
that $H^1(X,N_{\chi })\neq 0$. Define
$$\mathcal{V} :=\{x=(X,l,l',\mathcal{F}, \mathcal{F}'):(X,l,l')\in \mathcal{U},
\mbox{ }\mathcal{F} \mbox{ is a frame for } l,\mbox{ } \mathcal{F}' \mbox{ is a  frame
for }l'\}.$$ We may assume that for a generic $x\in \mathcal{U}$, the corresponding
pencils $l$ and $l'$ are base-point-free. Suppose that $\mathcal{U}$ has a component
of dimension $\alpha $. For any $x\in \mathcal{V}$,
$$T_{x}(\mathcal{V})\subseteq H^0(X, N_{\chi }),  \mbox{ and }\mbox { dim
}T_x(\mathcal{V})\geq \alpha +2\mbox{ dim }PGL(2)=\alpha +6.$$ If $\mathcal{K} _{\chi
}$ is the cuspidal sheaf of $\chi $ and $N_{\chi }'=N_{\chi }/\mathcal{K} _{\chi }$,
then according to \cite{AC1} Lemma 1.4, for a general point $x\in \mathcal{V}$ one has
that,
$$T_x(\mathcal{V})\cap H^0(X,\mathcal{K} _{\chi })=\emptyset ,$$
from which it follows that $\alpha \leq g-6.$ If not, one would have that 
$h^0(X,N_{\chi }')\geq g+1$, and therefore by Clifford's Theorem, $h^1(X,N_{\chi
})=h^1(X,N_{\chi }')=0$, which contradicts the definition of $\mathcal{U} $. Since
clearly $\mbox{dim }\mathcal{M}_{oct}> g-6$, we can assume that $h^1(C,N_{\phi })=0,$
for the general $[C]\in \mathcal{M}_{oct}$. Therefore, by taking cohomology in (11),
we get that $$H^1(C,N_{\phi _1 })=H^1(C,\mathcal{O}_C(2)),$$ where
$\mathcal{O}_C(1)=\phi _1^*\mathcal{O}_{\mathbb P^3}(1).$ By Serre duality,
\begin{equation}
H^1(C,\mathcal{O}_C(2))=0\Longleftrightarrow |K_C-2\mathfrak g^1_6-2 \mathfrak
h^1_6|=\emptyset.
\end{equation}
Since $K_C=5H-\sum_{i=1}^{6}(q_i'+q_i'')$, equation (12) becomes
\begin{equation}
|H+q_1'+q_1''+q_2'+q_2''-\sum_{i=3}^6(q_i'+q_i''|=\emptyset .
\end{equation}
If $L=\overline{p_1p_2}\subseteq \mathbb P^2$, we can write $\nu ^*(L)=
q_1'+q_1''+q_2'+q_2''+x+y+z+t$, and (13) is rewritten as
$$|2H-x-y-z-t-\sum_{i=3}^6(q_i'+q_i'')|=\emptyset.$$
So, one has to show that there are no conics passing through the nodes  $p_3,p_4,p_5$
and $p_6$ and also through the points in $L.\overline{C} -2p_1-2p_2$. Because
$[\overline{C}]\in U_{8,15}$ is general we may assume that $x,y,z$ and $t$ are
distinct, smooth points of $\overline{C} $. Indeed, if the divisor $x+y+z+t$ on
$\overline{C}$ does not consist of distinct points, or one of its points is a node, we
obtain that $\overline{C} $ has intersection number $8$ with the line $L$ at $5$
points or less. But according to \cite{DH}, the locus in the Severi variety
$$\{[X]\in U_{d,g}:X \mbox{ has total intersection number }m+3\mbox{ with a  line at
}m\mbox{ points }\}$$ is a divisor on $U_{d,g}$, so we may assume that
$[\overline{C}]$ lies outside this divisor. Now, if $x,y,z$ and $t$ are distinct and
smooth points of $\overline{C}$, a conic satisfying (13) would necessarily be a
degenerate one, and one gets a contradiction with the assumption that the nodes
$p_1,\ldots ,p_6$ of $\overline{C}$ are in general position.\hfill $\Box $
\newline
\textbf{Remark: } We have a nice geometric characterization of some of the strata 
$\mathcal{M}_i$. First, by using Zariski's Main Theorem for the birational projection
$\mathcal{G} ^1_6\rightarrow \mathcal{M} (1)$, one sees that $[C]\in \mathcal{M}
(1)_{sing}$ if and only if either $[C]\in \mathcal{M}(2)^0$, or $C$ possesses a
$\mathfrak g^1_6$ such that $\mbox{dim }|2\mathfrak g^1_6|\geq 3$. In the latter case,
the $\mathfrak g^1_6$ is a specialization of $2$ different $\mathfrak g^1_6$'s in some
family of curves, hence $\mathcal{M} (2)=\mathcal{M} (1)_{sing}$ (cf \cite{Co2}). As a
matter of fact, Coppens has proved that for $4\leq k\leq [(g+1)/2]$ and $8\leq g\leq
(k-1)^2$, there exists a $k$-gonal curve of genus $g$ carrying exactly $2$ linear
series $\mathfrak g^1_k$, so the general point of $\mathcal{M}  (2)$  corresponds to a
curve $C$ of genus $15$, having exactly $2$ base-point-free $\mathfrak g^1_6$'s.
Furthermore, using Coppens' classification of curves having many pencils computing the
gonality (see \cite{Co1}), we have that $\mathcal{M}(6)=\mathcal{M} _{oct}$, and
$\mathcal{M} (i)=\mathcal{M} _{sept}$, for each $i\geq 7$. 
\newline
\newline 
Now we are in a position to complete the proof of Theorem 1:\newline
\textsl{Proof of Theorem 1} According to (2), it will suffice to prove that there 
exists a smooth curve $[Y]\in \mathcal{M} _{23}$ which carries a $\mathfrak g^3_{20}$,
a $\mathfrak g^2_{17}$ but has no $\mathfrak g^1_{12}$'s. In the proof of Prop.5.4 we
constructed a stable curve of compact type  $[X]\in \overline{\mathcal{M}} _{23}$ such
that $[X]\in \overline{\mathcal{M} }^2_{17}\cap \overline{\mathcal{M} }^3_{20}$, but
$[X]\notin \overline{\mathcal{M}}^1_{12}$. If we prove that $[X]\in
\overline{\mathcal{M}^2_{17}\cap \mathcal{M} ^3_{20}}$, that is, there are smoothings
of $X$ which preserve both the $\mathfrak g^2_{17}$ and the $\mathfrak g^3_{20}$, we
are done. One can write  
$\overline{\mathcal{M} }^2_{17}\cap \overline{\mathcal {M}}^3_{20}=Y_1\cup  \ldots
\cup Y_s,$ where $Y_i$ are irreducible codimension $2$ subvarieties of
$\overline{\mathcal{M}} _{23}$.  Assume that $[X]\in Y_1.$ If $Y_1\cap \mathcal{M}
_{23}\neq \emptyset,$ then $[X]\in Y_1=\overline{Y_1\cap \mathcal{M} _{23}} \subseteq
\overline{\mathcal{M}^2_{17}\cap \mathcal{M} ^3_{20}}$, and the conclusion follows. So
we may assume that $Y_1\subseteq \overline{\mathcal{M} }_{23}-\mathcal{M} _{23}$.
Because $[X]\in \Delta _1-\bigcup_{j\neq 1}\Delta _j$, we must have $Y\subseteq \Delta
_1$. It follows that $\overline{\mathcal{M}} _{17}^2\cap \Delta_1$ and
$\overline{\mathcal{M}} _{20}^3\cap \Delta _1$ have $Y_1$ as a common component.
According to Prop.4.2, both intersections $\overline{\mathcal{M}} _{17}^2\cap \Delta
_1$ and $\overline{\mathcal{M} } _{20}^3\cap \Delta _1$ are irreducible, hence
$\overline{\mathcal{M}} _{17}^2\cap \Delta _1 =\overline{\mathcal{M}}_{20}^3\cap
\Delta _1=Y_1$, which contradicts Prop.5.3. Theorem 1 now follows.\hfill $\Box $
\section{The slope conjecture and $\mathcal{M}_{23}$}
In this final section we explain how the slope conjecture in the context of 
$\mathcal{M} _{23}$ implies that $\kappa (\mathcal{M} _{23})=2$, and then we present
evidence for this.\newline
\indent The slope of $\overline{\mathcal{M}}_g$ is defined as
$s_g:=\mbox{inf }\{a\in \mathbb R_{>0}:|a\lambda-\delta|\neq \emptyset\},\mbox{   
where  }\delta =\delta _0+\delta _1+\cdots +\delta _{[g/2]},\lambda \in
\mbox{Pic}(\overline{\mathcal{M}}_g)\otimes \mathbb R.$  Since $\lambda $ is big, it
follows that $s_g<\infty$. If $\mathbb E$ is the cone of effective divisors in
$\mbox{Div}(\overline{\mathcal{M}}_g)\otimes \mathbb R$, we define the slope function
$s:\mathbb E\rightarrow \mathbb R$ by the formula 
$$s_D:=\mbox{inf }\{a/b:a,b>0\mbox{ such that } \exists c_i\geq 0 \mbox{ with  }
[D]=a\lambda -b\delta -\sum_{i=0}^{[g/2]}c_i\delta _i\},$$
for an effective divisor $D$ on $\overline{\mathcal{M}}_g$. Clearly $s_g\leq s_D$  for
any $D\in \mathbb E$. When $g+1$ is composite, we obtain the estimate $s_{g}\leq
6+12/(g+1)$ by using the Brill-Noether divisors $\overline{\mathcal{M}}^r_{g,d}$, with
$\rho(g,r,d)=-1.$ 
\newtheorem{conj}{Conjecture}
\begin{conj}[\cite{HMo}]
We have that $s_g\geq 6+12/(g+1)$ for each $g\geq 3$, with equality when $g+1$ is
composite.
\end{conj}
\indent Harris and Morrison also stated (in a somewhat vague form) that for  composite
$g+1$, the Brill-Noether divisors not only minimize the slope among all effective
divisors, but they also single out those irreducible
$D\in \mathbb E$ with $s_D=s_g$.\newline
\indent The slope conjecture has been proved for $3\leq g\leq 11, g\neq 10$ (cf.\
\cite{HMo}, [CR3,4], \cite{Tan}). A strong form of the conjecture holds for $g=3$ and
$g=5$: on $\overline{\mathcal{M}}_3$ the only irreducible divisor of slope $s_3=9$ is
the hyperelliptic divisor, while on $\overline{\mathcal{M}}_5$ the only irreducible
divisor of slope $s_5=8$ is the trigonal divisor (cf. \cite{HMo}). Conjecture 1 would
imply that $\kappa (\mathcal{M}_g)=-\infty $ for all $g\leq 22$. For $g=23$, we
rewrite (1) as
\begin{equation}
nK_{\overline{\mathcal{M}}_{23}}=\frac{n}{c_{23,r,d}}[\overline{\mathcal{M}}^r_{g,d}]+
8n\,\delta _1+\sum _{i=2}^{11}\frac{(i(23-i)-4)}{2}n\,\delta _i\mbox{   }\mbox{ }
\mbox{  }(n\geq 1),
\end{equation}
(see Section 2 for the coefficients $c_{g,r,d}$). As Harris and Morrison suggest,  we
can ask the question whether the class of any $D\in \mathbb E$ with $s_D=s_g$ is
(modulo a sum of boundary components $\Delta _i$) proportional to
$[\overline{\mathcal{M}}^r_{23,d}]$, and whether the sections defining (multiples of)
$\overline{\mathcal{M}}_{23,d}^r$ form a maximal algebraically independent subset of
the canonical ring $R(\overline{\mathcal{M}}_{23})$. If so, it would mean that the
boundary divisor $8n\delta _1+(1/2)\sum _{i=2}^{11}n(i(23-i)-4)\delta _i$ is a fixed
part of $|nK_{\overline{\mathcal{M}}_{23}}|$. Moreover, using our independence result
for the three Brill-Noether divisors, it would follow that
$h^0(\overline{\mathcal{M}}_{23}, nK_{\overline{\mathcal{M}}_{23}})$ grows
quadratically in $n$, for $n$ sufficiently high and sufficiently divisible, hence
$\kappa (\mathcal{M}_{23})=2$. We would also have that $\Sigma \cap
\mathcal{M}_{23}=\mathcal{M}^1_{12}\cap \mathcal{M}^2_{17}\cap \mathcal{M}^3_{20}$,
with $\Sigma $ the common base locus of all the linear systems
$|nK_{\overline{\mathcal{M}}_{23}}|$.\newline
 \indent Evidence for these facts is of various sorts: first, one knows (cf. \cite{Ta}, \cite{CR3}) that $|nK_{\overline{\mathcal{M}}_{23}}|$ has a large fixed part in the boundary: for each $n\geq 1$, every divisor in $|nK_{\overline{\mathcal{M}}_{23}}|$ must contain $\Delta _i$ with multiplicity $16n$ when $i=1$, $19n$ when $i=2$, and $(21-i)n$ for $i=3,\ldots ,9$ or $11$. The results for $\Delta _1$ and $\Delta _2$ are optimal since these multiplicities coincide with those in (14). Note that $[\Delta _1]=2\delta _1$.\newline \indent Next, one can show that certain geometric loci in $\mathcal{M}_{23}$ which are contained in all three Brill-Noether divisors, are contained in $\Sigma $ as well. The method is based on the trivial observation that for a family $f:X\rightarrow B$ of stable curves of genus $23$ with smooth general member, if $B.K_{\overline{\mathcal{M}}_{23}}<0$ (or equivalently, slope$(X/B)=\delta _B/\lambda _B>13/2)$, then $\phi (B)\subseteq \Sigma$, where
$\phi:B\rightarrow \overline{\mathcal{M}}_{23}, \phi(b)=[X_b]$, is the associated 
moduli map. We have that:
\newline
$\bullet$ One can fill up the $d$-gonal locus $\overline{\mathcal{M}}^1_d$ with 
families $f:X\rightarrow B$ of stable curves of genus $g$ such that slope$(X/B)$ is
$8+4/g$ in the hyperelliptic case, $36(g+1)/(5g+1)$ in the trigonal case, and
$4(5g+7)/(3g+1)$ in the tetragonal case (cf. \cite{Sta}). For $g=23$ it follows that
$\mathcal{M}^1_4\subseteq \Sigma $. Note that this result is not optimal if we believe
the slope conjecture since we know that $\mathcal{M}^1_8\subseteq
\mathcal{M}^1_{12}\cap \mathcal{M}^2_{17}\cap \mathcal{M}^3_{20}$. (The inclusion
$\mathcal{M}^1_8\subseteq \mathcal{M}^3_{20}$ is a particular case of a result from
\cite{CM}.)\newline
$\bullet$ We take a pencil of nodal plane curves of degree $d$ with $f$ assigned
nodes in general position such that ${d-1\choose 2}-f=23$, and with $b$ base points, 
where $4f+b=d^2$. After blowing-up the base points, we have a pencil
$Y\rightarrow \mathbb P^1$ with fibre $[Y_t]\in \overline{\mathcal{M}}^2_d$.
For this pencil $\lambda _{\mathbb P^1}= \chi(\mathcal{O}_Y)+23-1=23$ and
$\delta_{\mathbb P^1}=c_2(Y)+88=91+b+f$. The condition $\delta_{\mathbb P^1}/\lambda
_{\mathbb P^1}>13/2$ is satisfied precisely when $d\leq 10$, hence taking into account
that such pencils fill up $\mathcal{M}^2_d$, we obtain that
$\mathcal{M}^2_{10}\subseteq \Sigma$. Note that $\mathcal{M}^2_{10}\subseteq
\mathcal{M}^1_8$, and as mentioned above, the $8$-gonal locus is contained in the
intersection of the Brill-Noether divisors.\newline
$\bullet$ In a similar fashion we can prove that $\mathcal{M}_{23,\gamma }(2)$,  the
locus of curves of genus $23$ which are double coverings of curves of genus $\gamma $
is contained in $\Sigma $ for $\gamma \leq 5$.
\newline
\indent The fact that the slopes of other divisors on 
$\overline{\mathcal{M}}_{23}$ (or on $\overline{\mathcal{M}}_g$ for arbitrary $g$)
consisting of curves with special geometric characterization, are larger than
$6+12/(g+1)$, lends further support to the slope hypothesis. In another paper we will
compute the class of various divisors on $\overline{\mathcal{M}}_{23}$: the closure in
$\overline{\mathcal{M}}_{23}$ of the locus
$$
D_e:=\{[C]\in \mathcal{M}_{23}:\mbox{there exists }  l\in G^1_{(23+e)/2}(C)\mbox{
and } p\in C,\mbox{ such that } w^l(p)\geq e\},$$ for $3\leq e\leq 19$ and $e$ odd,
and the closure of the locus
$$
\{[C]\in \mathcal{M}_{23}:C\mbox{ has a }\mathfrak g^2_{18} \mbox{ with a
}5\mbox{-fold point, i.e. }\exists D\in C^{(5)}\mbox{ such that }\mathfrak
g^2_{18}(-D)=\mathfrak g^1_{13}\}.$$ In each case we will show that the 
slope estimate
holds. \vskip10pt
\noindent Korteweg-de Vries Institute for Mathematics, University of Amsterdam\newline
Plantage Muidergracht 24, 1018 TV Amsterdam, The Netherlands\newline
e-mail: {\tt farkas@wins.uva.nl }

\end{document}